\documentclass[11pt,english,oneside]{amsart}
\usepackage[T1]{fontenc}
\usepackage[latin9]{inputenc}
\usepackage{geometry}
\usepackage{amssymb}
\usepackage{color}

\usepackage{graphicx,color}
\usepackage{amsmath, amssymb, graphics}

\usepackage{graphicx}

\makeatletter
\numberwithin{equation}{section} 
\numberwithin{figure}{section} 
  \@ifundefined{theoremstyle}{\usepackage{amsthm}}{}
  \theoremstyle{plain}
  \newtheorem{thm}{Theorem}[section]
  \theoremstyle{plain}
  \newtheorem{cor}[thm]{Corollary}
  \theoremstyle{plain}
  
  \theoremstyle{Remark}
  \newtheorem{rem}[thm]{Remark}
  \theoremstyle{remark}
  
  \theoremstyle{plain}
  \newtheorem{lem}[thm]{Lemma}

\usepackage{geometry}

\def\bfR#1{{\bf R}^#1}

\def\com#1{ \hbox{#1}}

\smallskip
\def\<{{\langle }}
\def\>{{\rangle }}

\makeatother

\usepackage{babel}

\def\bfR#1{{\bf R}^#1}

\def\com#1{ \quad\hbox{#1}\quad}

\smallskip
\def\<{{\langle }}
\def\>{{\rangle }}

\makeatother

\begin{document}

\title{Embedded CMC hypersurfaces on hyperbolic spaces}

\author{ Oscar M. Perdomo }

\date{\today}

\curraddr{Department of Mathematics\\
Central Connecticut State University\\
New Britain, CT 06050\\
}

\email{ perdomoosm@ccsu.edu}

\begin{abstract}

In this paper we will prove that for every integer $n>1$, there exists a real number $H_0<-1$ such that  every  $H\in (-\infty,H_0)$ can be realized as the mean curvature of a embedding of  $H^{n-1}\times S^1$  in the $n+1$-dimensional spaces $H^{n+1}$. For $n=2$ we explicitly compute the value $H_0$. For a general value $n$, we provide function $\xi_n$ defined on $(-\infty,-1)$, which is easy to compute numerically, such that, if $\xi_n(H)>-2\pi$, then, $H$ can be realized as the mean curvature of a embedding of  $H^{n-1}\times S^1$  in the $n+1$-dimensional spaces $H^{n+1}$.
\end{abstract}

\subjclass[2000]{58E12, 58E20, 53C42, 53C43}

\maketitle

\section{Introduction and preliminaries}

Here we will be considering the following model of the hyperbolic space,

$$H^{n+1}=\{\, x\in \bfR{{n+2}}: x_1^2+\dots+x_{n+1}^2-x_{n+2}^2=-1\, \} $$

where the space $\bfR{{n+2}}$ is endowed with the following inner product

$$\<v,w\>=v_1w_1+\dots + v_{n+1}w_{n+1}-v_{n+2}w_{n+2}  \com{for $v=(v_1,\dots,v_{n+2})$ and $w=(w_1,\dots,w_{n+2})$}  $$

In \cite{P} we proved the following theorem that shows that $S^{n-1}\times {\bf R}$ can be embedded in the hyperbolic space with constant mean curvature.

\begin{thm}

Let $g_{C,H}:{\bf R}\to {\bf R}$ be a positive solution of the equation

\begin{eqnarray}\label{eq hyperbolic}
 (g^\prime)^2+g^{2-2n}+(H^2-1)g^2+2Hg^{2-n}=C
 \end{eqnarray}

associated with a non negative $H$  and a positive constant $C$. If $\mu,\lambda,r,\theta:{\bf R}\to {\bf R}$  are defined by

$$r=\frac{g_{C,H}}{\sqrt{C}},\quad  \lambda=H+g_{C,H}^{-n},\, \mu=nH-(n-1)\lambda=H-(n-1)g_{C,H}^{-n} \com{and}\theta(u)=\int_0^u\frac{r(s)\lambda(s)}{1+r^2(s)}ds$$

then, the map $\phi:S^{n-1}\times {\bf R}\to H^{n+1}$ given by

\begin{eqnarray}\label{the immersions hyperbolic}
\phi(y,u)=(\, r(u)\, y,\sqrt{1+r(u)^2}\, \sinh(\theta(u)),\sqrt{1+r(u)^2} \, \cosh(\theta(u))\, )
\end{eqnarray}

defines an embedded hypersurface in $H^{n+1}$ with constant mean curvature $H$. Moreover, if
$H^2>1$,  the embedded manifold defined by $(\ref{the immersions hyperbolic})$ admits the group $O(n)\times Z$ in its group of isometries, where $Z$ is the group of integers.
\end{thm}

The existence of the previous examples just as immersions were studied in \cite{SI} as Delaunay-type hypersurfaces of the hyperbolic space and also in \cite{D-D} as rotational hypersurfaces of spherical type. In this paper we will prove that a subfamily of the family of immersions named as {\it rotational hypersurfaces of hyperbolic type} in  \cite{D-D}  provides different ways to embed the manifold $H^{n-1}\times S^1$ in the $n+1$-dimensional hyperbolic space.

\section{Embedded hyperbolic type rotational surfaces in $H^3$ }

It is not difficult to show that the function

$$\xi:(-\infty,-1)\to {\bf R}\com{given by}
\xi(H)=\int_0^\pi\frac{\sqrt{2}\, H\, dt}{\sqrt{2H^2+\sin(2t)-1}}$$

is decreasing, $\lim_{H\to-\infty}\xi(H)=-\pi$ and $\xi(H)<-2\pi$ for values of $H$ close to $-1$. The previous observations guarantee the existence of a unique $H_0$ such that $\xi(H_0)=-2\pi$. A numerical computation shows that,

$$H_0\simeq -1.0158136657178574$$

In this section we will show that every $H<H_0$ can be realized as the mean curvature of a hyperbolic type rotational embedded constant mean curvature surface in the hyperbolic three dimensional space. Let us state and prove the main and only theorem in this section.

\begin{thm}
For any $H<-1$ and $C\in \, (C_1,0\, )$ where $C_1=\, 2 (H + \sqrt{-1 + H^2})$, let us define $f:{\bf R}\to {\bf R}$ by

$$f(t)= \sqrt{\frac{C - 2 H +
 \sqrt{4 + C^2 - 4 C H} \, \sin(\,2 \sqrt{H^2-1\, }\, t\, )}{2H^2-2}}
$$

If we define,

$$r(t)=\frac{f(t)}{\sqrt{-C}}\com{and} \lambda(t)=H+(f(t))^{-2}$$

then, the function  $\frac{\lambda(t)\, r(t)}{r^2(t)-1}$ is a smooth function everywhere and if we define

$$\theta(t)=\int_0^t\frac{\lambda(s)\, r(s)}{r^2(s)-1}\, ds$$

then, the map

\begin{eqnarray}\label{the immersions hyperbolic 2 case n=2}
\phi(y,u)=(\, \sqrt{r(u)^2-1}\, \cos(\theta(u)),\, \sqrt{r(u)^2-1} \, \sin(\theta(u)),\,  r(u) \, \sinh(v), r(u) \, \cosh(v) )
\end{eqnarray}

defines and immersion from $\bfR{2}$ to $H^3$. We also have that for every $H<-1$ there exist infinitely many choices of $C$ such that the immersion $\phi$ is periodic in the variable $u$ and therefore it defines immersions from ${\bf R}\times S^1$ to $H^3$. Moreover, we have that for every $H<H_0$, there exists a value $C$ such that $\phi$ defines an embedding from ${\bf R}\times S^1$ to $H^3$.
\end{thm}

\begin{proof}
Since $H<-1$ and $C\in \, (\,C_1 ,0\, )$, we have that the function $f$ is a real-value $T$-periodic function, that oscillates from $t_1$ to $t_2$ where,

$$
t_1,\, t_2=\sqrt{\frac{C - 2 H \pm
 \sqrt{4 + C^2 - 4 C H} \, }{2H^2-2}}\com{and} T=\frac{\pi}{\sqrt{H^2-1}}
$$

A direct computation shows that

$$
(f^\prime)^2+f^{-2}+(H^2-1)f^2+2H=C
$$

The equation above shows that the function $r(t)$ satisfies the following identity

\begin{eqnarray}
(r^\prime)^2+\lambda^2\, r^2= r^2-1
\end{eqnarray}

This equation shows that $r(t)\ge 1$, moreover, it shows that $r(t^\star)=1$, if and only if $\lambda(t^\star)=0$ and $r^\prime(t^\star)=0$. These last two conditions imply that $t^\star$ is a zero with multiplicity $2$ of the function $\lambda$. We can easily see that $t^\star$ is a zero  with multiplicity $2$ of the function $r^2-1$. Since the function $r$ is analytic, we get that the function $\frac{\lambda(s)\, r(s)}{r^2(s)-1}$ is smooth near $t^\star$, therefore it is smooth everywhere.  A direct computation shows that

\begin{eqnarray*}
\frac{\partial \phi}{\partial u}= \frac{r\, r^\prime}{\sqrt{r^2-1}} (\cos(\theta),\sin(\theta),0,0)+\frac{r\, \lambda}{\sqrt{r^2-1}}\, (-\sin(\theta),\cos(\theta),0,0)+ r^\prime(0,0,\, \sinh(v),\cosh(v))
\end{eqnarray*}

 and

 $$\frac{\partial \phi}{\partial v}=r(u)\, (0,0,\, \cosh(v),\sinh(v))$$

 It is not difficult to prove that the map

$$ \nu=-r\lambda\,   (0,0,\sinh(v),\cosh(v))-\frac{r^2\, \lambda}{\sqrt{r^2-1}} \, (\cos(\theta),\sin(\theta),0,0\, ) + \frac{r^\prime}{\sqrt{r^2-1}} \, (-\sin(\theta),\cos(\theta),0,0\, )$$

is a Gauss map of the immersion $\phi$. It follows that the immersion $\phi$ has constant mean curvature $H$ by noticing that

$$\frac{\partial \nu}{\partial v}=-\lambda\frac{\partial \phi}{\partial v}\com{and}
\frac{\partial \nu}{\partial u}=-(2 H-\lambda)\frac{\partial \phi}{\partial u}$$

Let us define the function $K$ that depends on $H$ and $C$, by

$$K(C,H)=\int_0^T \frac{\lambda(s)\, r(s)}{r^2(s)-1}\, ds$$

A direct computation shows that for every fixed $H$ we have,

\begin{eqnarray}\label{limit n=2}
\lim_{C\to C_1}\, K(C,H)=-\, \pi\, \sqrt{2 - \frac{2 H}{\sqrt{H^2-1}}} =b_2(H)\com{and}\lim_{C\to 0}\, K(C,H)=0
\end{eqnarray}

Since $H<-1$ we have that $b_2(H)<-2\pi$. Using the limits in (\ref{limit n=2}) we get that for any fixed value $H<-1$ and for every positive integer $m$, there exists a real number $C^\star$ between $C_1$ and $0$ such that $K(C^\star,H)=-\frac{2\pi}{m}$. Since the function $\theta$ satisfies that

$$ \com{For any integer $j$ and $u\in [jT,(j+1)T]$ we have that }\theta(u)=jK+\theta(u-jT) $$

we get that if we choose the value $C^\star$, we get that  $\theta(m T)= -2\pi$ and therefore the immersion $\phi(u,v)$ will be $mT$-periodic in the variable $u$ and it will define an immersion from ${\bf R}\times S^1$ to $H^3$. Let us prove that for every $H<H_0$ there exists an embedding from ${\bf R}\times S^1$ to $H^3$.
By using the definition of the function $\lambda$ and the expression for the bounds $t_1$ and $t_2$ of the function $f$, we have that for a given $H$, the function $\lambda<0$ if and only if $C_1<C<\frac{1}{H}$. Notice that if $\lambda$ is always negative,  then the function $\theta$ is strictly decreasing, and in particular it is one to one. A direct computation shows that

$$K(\frac{1}{H},H)=\int_0^T\frac{H\, \sqrt{2H^2-2}}{\sqrt{2 H^2-1+\sin(\, 2\sqrt{H^2-1}\, \, s)\, }}\, ds= \int_0^\pi\frac{H\, \sqrt{2}}{\sqrt{2 H^2-1+\sin(\, 2 t)}\, }$$

As pointed out at the beginning of this section, the function $\xi(H)=K(\frac{1}{H},H)$ is decreasing and the limit when $H\to -\infty$ is $-\pi$. Therefore for any $H<H_0$ there exists a $C^\star$ between $C_1$ and $\frac{1}{H}$ such that $K(C^\star,H)=-2\pi$. By the way we picked $C^\star$ we get that the function $\theta$ is strictly decreasing and $\theta(T)=-2\pi$, these two conditions guarantee that the immersion $\phi(u,v)$ is $T$-periodic and injective in ${\bf R}\times (0,T)$, therefore $\phi$ defines an embedding from ${\bf R}\times S^1$ to $H^3$. This completes the proof of the theorem.
\end{proof}

\vskip.5cm

\subsection{Graph of some profile curves}

The examples described above are obtained by doing a hyperbolic rotation of the profile curve

$$\alpha(t)= (\sqrt{r^2(t)-1}\, \cos(\theta(t)), \, \sqrt{r^2(t)-1}\,\sin(\theta(t))\, )$$

We will show the graphs of a profile curve that corresponds to an embedded example and two profile curves corresponding to immersed examples, all of them represent examples with constant mean curvature $H=-1.1$. To finish the section we will show one of the numerical difficulties to do the graph. This difficulty is the fact that the angle function $\theta$ moves  a lot in a small variation of the parameter $t$, during this small variation of parameter $t$, the radius function $\sqrt{r^2(t)-1}$ is very close to zero. We will show this fact by graphing the function $\theta^\prime(t)$, first by limiting the codomain to some values close to zero, and then by showing the whole graph of $\theta^\prime$.

\vfill
\eject

\begin{figure}[h]
\centerline{\includegraphics{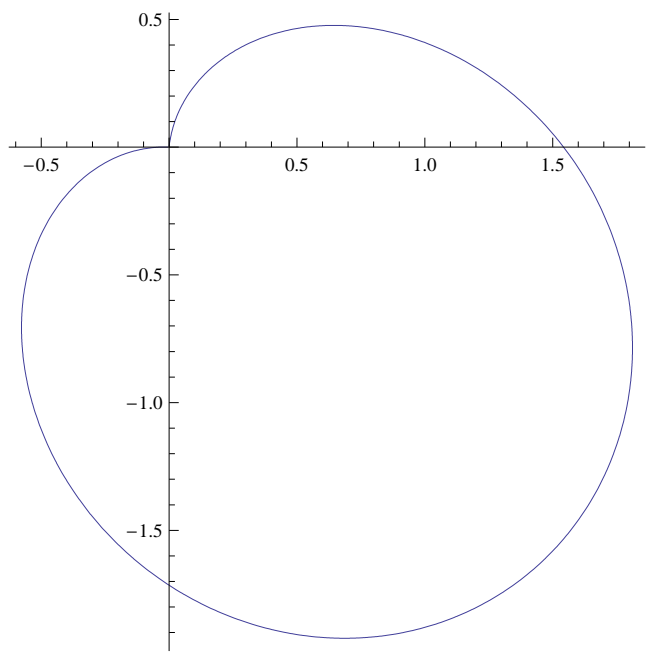}}
\caption{Profile curve for a surface with CMC $H=-1.1$, in this case the surface is embedded and $C= -0.9091743461769703$ and $K=-2\pi$ }
\end{figure}

\begin{figure}[h]
\centerline{\includegraphics{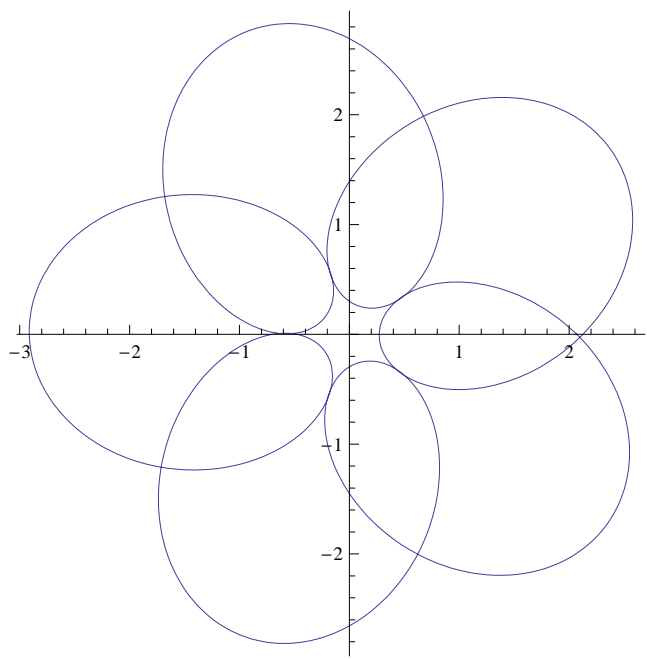}}
\caption{Profile curve for a surface with CMC $H=-1.1$, in this case $C= -0.6835660909345689$ and $K=-\frac{2\pi}{5}$ }
\end{figure}

\vfill
\eject

\begin{figure}[h]
\centerline{\includegraphics{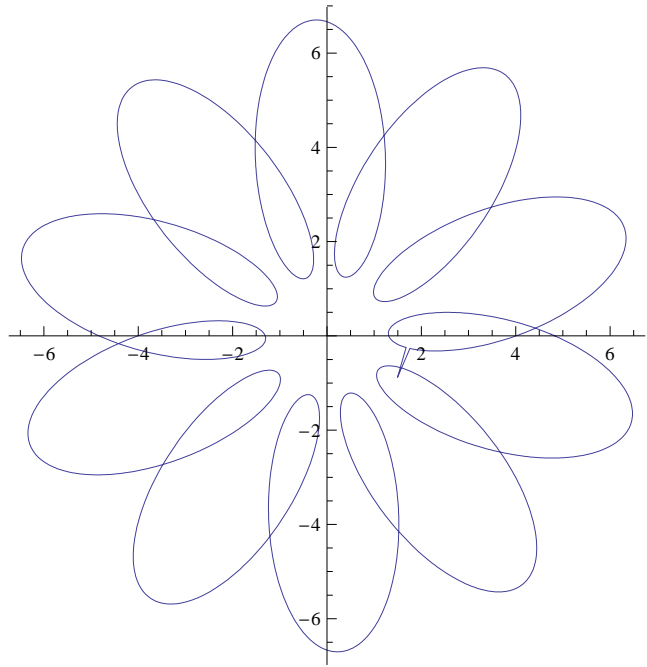}}
\caption{Profile curve for a surface with CMC $H=-1.1$, in this case $C= -0.19607165524075582$ and $K=-\frac{2\pi}{10}$ }
\end{figure}

\vfill
\eject

\begin{figure}[h]
\centerline{\includegraphics{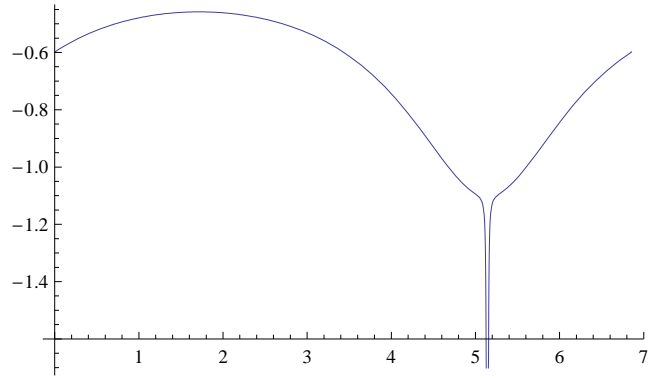}}
\caption{Graph of the function $\theta^\prime$ associated with the embedded example which profile curve is shown above, in this case just part of the graph is shown}
\end{figure}

\vskip1cm

\begin{figure}[h]
\centerline{\includegraphics{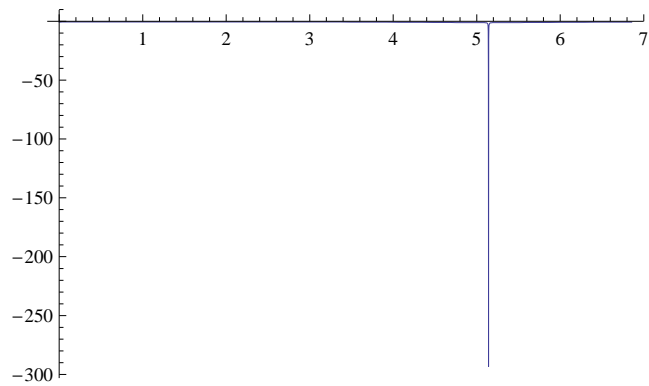}}
\caption{Graph of the function $\theta^\prime$ example associated with the embedded which profile curve is shown above.}
\end{figure}

\vfill
\eject

\section{Embedded solutions in hyperbolic spaces.}

It is well known that the existence of CMC hypersurfaces in hyperbolic spaces relies on the existence of solutions of the following differential equation,

\begin{eqnarray*}
 (g^\prime)^2+g^{2-2n}+(H^2-1)g^2+2Hg^{2-n}=C
\end{eqnarray*}

It is not difficult to check that, when $H<-1$,  it is possible to obtain solutions of this equation associated with negative values of  $C$. This CMC examples produced by these solutions when $C<0$ correspond to those  named as rotational hyperbolic type in \cite{D-D}. Similar arguments as those shown in \cite{P} will give us explicit immersions for such a choice of the constant $C$. The following sequences of statements  tell us how to pick the negative values of $C$ to obtain solutions in the case that $H<-1$ and several other properties that will be useful in the proof of the main theorem in this paper.

\begin{rem}\label{properties of q}
The function $q:(0,\infty)\to {\bf R}$  defined by $q(v)=C - v^{2 - 2 n} + (1 - H^2) v^2 - 2 H v^{2 - n}$, where $H<-1$ and $C<0$, has the following properties:
\begin{enumerate}

\item

The positive real number $v_0$ given by

\begin{eqnarray*}\label{value of v0}
v_0= (\, \frac{H(n-2) + \sqrt{4 - 4 n + H^2 n^2}}{2 H^2-2}\, )^\frac{1}{n}= (\, \frac{2(n-1)}{\sqrt{4 - 4 n + H^2 n^2}-H(n-2)}\, )^\frac{1}{n}
\end{eqnarray*}

is the only positive critical point of $q$.

\vskip.5cm
\item

$p(v)=v^{2n-2}q(v)$ is a polynomial with even degree, negative leading coefficient and $p(0)=-1$.

\vskip.5cm

\item
Since $q^{\prime}(v)>0$ if $v<v_0$, $q^{\prime}(v)<0$ if $v>v_0$, and $q(v_0)=C-C_0$  where

\begin{eqnarray}\label{value of c0}
C_0=  n\, \frac{H^2 n-2 + H \sqrt{4 - 4 n + H^2 n^2}}{(\, H (n-2) + \sqrt{4 - 4 n + H^2 n^2}\, )^{\frac{2n-2}{n}}} \, (2 H^2-2)^{\frac{n-2}{n}}
\end{eqnarray}

then, $q$ has exactly $2$ roots whenever $0>C>C_0$.

\vskip.5cm

\item

The functions $t_1,t_2:(C_0,0)\times(-\infty,-1)\to (0,\infty)$ defined by the equations

\begin{eqnarray}\label{definition of ti}
q(t_1(C,H))=0\quad q(t_2(C,H))=0\com{with} t_1(C,H)<t_2(C,H)
\end{eqnarray}

are smooth, $t_1(C,H)$ is decreasing with respect to $C$, $t_2(C,H)$ is increasing with respect to $C$ and
the limit of both functions when $C\to C_0$ is $v_0$.

\vskip.5cm

\item

Since the roots of $q$ when $C=0$ are $v_1=\frac{1}{(1 - h)^\frac{1}{n}}$ and $v_2=\frac{1}{(-1 - h)^\frac{1}{n}}$
then for any fixed $H$ the derivative of the functions $t_1$ and $t_2$ defined on $(C_0,0)$ never vanish and

$$\lim_{C\to 0}t_1(C)=v_1,\quad \lim_{C\to 0}t_2(C)=v_2\com{and}\lim_{C\to C_0}t_1(C)=\lim_{C\to C_0}t_2(C)=v_0$$

\vskip.5cm
\item

The following identities are true,

$$  \lambda_1=H+v_0^{-n}= \frac{n\,  H  + \sqrt{ H^2 n^2-4(n-1)}\, }{2 (n-1)}\, <\, 0 \com{and} \lambda_2=H+v_1^{-n}=1 $$

\vskip.5cm

\item

For a fixed $H<-1$, the previous two items guarantee the existence of a unique $\tilde{C}(H)\in(C_0,0)$ such that $t_1^{\star}(H)=t_1(\tilde{C}(H),H)$ satisfies that

$$H+(t_1^{\star}(H))^{-n}=0$$

The equality above defines a smooth function $\tilde{C}:(-\infty,-1)\to (C_0,0)$

\vskip.5cm

\item

We can explicitly compute the function $\tilde{C}$ by noticing first that for that special value of $C$, the number $t_1=(-H)^{\frac{-1}{n}}$ must be a root of the function $q$, therefore $q((-H)^{\frac{-1}{n}})$ must be zero, i.e,

$$ q((-H)^{\frac{-1}{n}})=C+(-H)^{-\frac{2}{n}  }=0$$

Therefore, $\tilde{C}(H)=-(-H)^{-\frac{2}{n}}$

\vskip.5cm

\item

The function $\tilde{q}(v)=-\frac{1}{C}q(\sqrt{-C} v)$ has the following expression

$$\tilde{q}(v)=-1 - (-C)^{-n} v^{2 - 2 n} + v^2 (1 - H^2 - 2 H (\sqrt{-C} v)^{-n})$$

Moreover, by the definition of $\tilde{q}$ and the properties of the function $q$ we have that $\tilde{q}$, for any $C\in(C_0,C)$, the only 2 positive roots of $\tilde{q}$  are

 $$\tilde{t_1}(C,H) = \frac{t_1(C,H)}{\sqrt{-C}}        \com{and}      \tilde{t_2}(C,H) =\frac{t_2(C,H)}{\sqrt{-C}}    $$

Therefore, we have that $\tilde{t_1}(\tilde{C},H)=\frac{(-H)^{-\frac{1}{n}}}{\sqrt{(-H)^{-\frac{2}{n}}}}=1$

\vskip.5cm
\item

A direct computation shows that when $C=\tilde{C}$, the polynomial $\tilde{q}$, reduces to the polynomial $Q$ given by,

$$Q=-1 + v^2 - H^2 v^2 - H^2 v^{2 - 2 n} + 2 H^2 v^{2 - n} $$

It is not difficult to check that, when $n>2$, for any positive $\epsilon$, $\lim_{H\to -\infty}Q(1+\epsilon)=-\infty$, therefore we have that

\begin{eqnarray}\label{values of tilde ti}
\tilde{t_1}(\tilde{C},H)=1\com{and} \lim_{H\to-\infty}\tilde{t_2}(\tilde{C},H)=1
\end{eqnarray}

\vskip.5cm
\item

Let us define the function $h:(0,\infty)\to {\bf R}$ by

$$ h(v)=\frac{2 H v^{1 - n} (-1 + v^n)}{v^2-1}=2 H\, v^{1 - n}\, \frac{1+v+\dots+v^{n-1}}{1+v}$$

and the function $\xi_n:(-\infty,-1)\to {\bf R}$

$$\xi_n(H)=\int_1^{\tilde{t_2}(\tilde{C},H)} \frac{h_n(v)}{\sqrt{Q(v)}}dv$$

\vskip.5cm
\item

A direct computation shows that

\begin{eqnarray}
\tilde{a}=-\frac{1}{2} Q^{\prime\prime}(1)=n^2H^2-1
\end{eqnarray}

Therefore using a small modification of lemma 5.1 and its corollary in \cite{P} we get that

\begin{eqnarray}\label{minuspilimit}
\lim_{H\to-\infty} \xi_n(H)=-\pi
\end{eqnarray}

\end{enumerate}

\end{rem}

\vskip.7cm

\begin{thm}\label{case n immersion}

Let $g:{\bf R}\to {\bf R}$ be a positive solution of the equation

\begin{eqnarray}\label{eq hyperbolic}
 (g^\prime)^2+g^{2-2n}+(H^2-1)g^2+2Hg^{2-n}=C
 \end{eqnarray}

associated with a negative constant $C$. If $\mu,\lambda,r,\theta:{\bf R}\to {\bf R}$  are defined by

$$r=\frac{g}{\sqrt{-C}},\quad  \lambda=H+g^{-n},\, \mu=nH-(n-1)\lambda=H-(n-1)g^{-n} \com{and}\theta(u)=\int_0^u\frac{r(s)\lambda(s)}{r^2(s)-1}ds$$

then, the map $\phi_{C,H}:H^{n-1}\times {\bf R}\to H^{n+1}$ given by

\begin{eqnarray}\label{the immersions hyperbolic 2}
\phi_{C,H}(y,u)\, =\,(\, \sqrt{r(u)^2-1}\, \cos(\theta(u)),\, \sqrt{r(u)^2-1} \, \sin(\theta(u)),\,  r(u) \, y )
\end{eqnarray}

defines an immersed hypersurface in $H^{n+1}$ with constant mean curvature $H$. We also have that when $H<-1$, the function $g$ is periodic and if we denote its period by $T$, then, $\phi_{C,H}$ defines an immersion from $H^{n-1}\times S^1$ to $H^n$ whenever

\begin{eqnarray}\label{expression for K}
K(C,H)=\int_0^T\frac{r(s)\lambda(s)}{r^2(s)-1}ds=-\frac{2 k \pi}{m}
\com{for some pair of positive integers $k$ and $m$}
\end{eqnarray}

Moreover, we have that anytime $\xi_n(H_1)>-2 \pi$, where $\xi_n$ is the function defined in item (11) in Remark (\ref{properties of q}), then, there exists a constant $C$ such that the immersion $\phi_{C,H_1}$ defines an embedding from $H^{n-1}\times  S^1$ to $H^{n+1}$.

\end{thm}

\begin{proof}

A direct computation shows the following identities,

$$ (r^\prime)^2+r^2 \lambda^2=r^2-1,\com{and} \lambda r^\prime +r\lambda^\prime=\mu r^\prime$$

Let us define

$$B_2(u)=(\cos(\theta(u)),\sin(\theta(u)),0,\dots,0)\com{and} B_3(u)=(-\sin(\theta(u)),\cos(\theta(u)),0,\dots,0)$$

Notice that $\<B_2,B_2\>=1$,  $\<B_3,B_3\>=1$, $\<B_2,B_3\>=0$,  $B_2^\prime=\frac{r\lambda}{r^2-1}B_3$ and
$B_3^\prime=-\frac{r\lambda}{r^2-1}B_2$, moreover, we have that the map $\phi=\phi_{C,H}$ can be written as

$$\phi=\, r\, (0,0,y)+\sqrt{r^2-1} \, B_2$$

A direct verification shows that $\<\phi,\phi\>=-1$ and that

$$\frac{\partial{\phi}}{\partial u}=r^\prime\, (0,0,y)+\frac{rr^\prime}{\sqrt{r^2-1}} \, B_2+
\frac{r \lambda}{\sqrt{r^2-1}} \, B_3$$

is a unit vector, i.e, $\<\frac{\partial{\phi}}{\partial u},\frac{\partial{\phi}}{\partial u}\>=1$.
We have that the tangent space of the immersion at $(y,u)$ is given by

$$T_{\phi(y,u)}=\{(v,0,0)+s\, \frac{\partial \phi}{\partial u}: \<v,y\>=0\com{and} s\in{\bf R}\}$$

A direct verification shows that the map

$$\nu=-r\lambda\,   (0,0,y)-\frac{r^2\, \lambda}{\sqrt{r^2-1}} \, B_2 + \frac{r^\prime}{\sqrt{r^2-1}} \, B_3$$

satisfies that $\<\nu,\nu\>=1$, $\<\nu,\frac{\partial \phi}{\partial u}\>=0$ and, for any $v\in\bfR{n}$ with $\<v,y\>=0$, we have that $\<\nu,(v,0,0)\>=0$. It then follows that $\nu$ is a Gauss map of the immersion $\phi$. The fact that the immersion $\phi$ has constant mean curvature $H$ follows because, for any unit vector $v$ in $\bfR{n}$ perpendicular to $y$, we have that

$$\beta(t)=(\, 0,\, 0\, ,\,  r\cosh(t)\, y+r\sinh(t)\, v \,)+\sqrt{r^2-1}\, B_2=\phi(\cosh(t)y+\sinh(t)v,u)$$

satisfies that $\beta(0)=\phi(y,u)$, $\beta^\prime(0)=rv$ and

$$\frac{d \nu(\beta(t))}{dt}\big{|}_{t=0} = d\nu(rv)=-r\lambda\, v$$

Therefore,  $\lambda$ is a principal curvature with multiplicity $n-1$. Now, since $\<\frac{\partial \nu}{\partial u},(v,0,0)\>=0$ for every $(v,0,0)\in T_{\phi(y,u)}$, we have that $\frac{\partial \phi}{\partial u}$ defines a principal direction, i.e. we must have that $\frac{\partial \phi}{\partial u}$ is a multiple of $\frac{\partial \phi}{\partial u}$.  A direct verification shows that,

$$\<\frac{\partial \nu}{\partial u},y\>=-\lambda^\prime \, r-\lambda r^\prime=-\mu\, r^\prime=-(nH-(n-1)\lambda)r^\prime$$

We also have that $\<\frac{\partial \phi}{\partial u},y\>=r^\prime$, therefore,

$$\frac{\partial \nu}{\partial u}= d\nu(\frac{\partial \phi}{\partial u})=-\mu\, \frac{\partial \phi}{\partial u}=-(nH-(n-1)\lambda) \frac{\partial \phi}{\partial u}$$

It follows that the other principal curvature is $nH-(n-1)\lambda$. Therefore $\phi$ defines an immersion with constant mean curvature $H$, this proves the first item in the Theorem. The fact that the map \label{the immersions hyperbolic 2} defines an immersion from $H^{n-1}\times S^1$ whenever $K(C,H)=-\frac{2 k \pi}{m}$, follows from the following property

$$ \com{For any integer $j$ and $u\in [jT,(j+1)T]$ we have that }\theta(u)=jK+\theta(u-jT) $$

which implies that the map $\phi$ is periodic in the variable $u$, with period $m T$. Let us prove the embedding part of the theorem. In this part of the proof we will be using the functions and constants

$$q,\, \tilde{q},\, Q,\, \xi_n,\,  t_1,\, t_2,\, \tilde{t_1},\, \tilde{t_2},\, \tilde{C},\, C_0,\, \com{and} v_0$$

defined in Remark (\ref{properties of q}). Let us start by noticing that the differential equations for the functions $g$ and $r$ can be written as

$$(g^\prime)= q(g) \com{and} (r^\prime)^2=\tilde{q} (r)$$

It follows that, in order to obtain a solution $g$ of this differential equation, we need that $C>C_0$ and, once we have the solution $g$ associated with the number $C$ and $H$, this solution $g$ varies from $t_1(C,H)$ to $t_2(C,H)$. Since we know the maximum and the minimum of the function $g$ in terms of $C$ and $H$, we can verify that anytime $C<\tilde{C}=-(-H)^{-\frac{2}{n}}$, the function $\lambda$ is negative, we also have that when $C=\tilde{C}$, $0$ is the maximum of the function $\lambda$. The previous affirmation guarantees that anytime $C\in(C_0,\tilde{C})$, the function $\theta$ is one to one. By doing the substitution $v=g(s)$ in the integral  $K(C,H)$, we get that

$$K(C,H)=\int_{t_1(C,H)}^{t_2(C,H)}  \frac{2 \sqrt{-C} \,  (1 + H v^n)\, v^{1 - n}}{(C + v^2)\, \sqrt{q(v)}}\, dv$$

In the previous expression we have used the symmetry of the function $g$, and therefore the symmetries of the functions $r$ and $\lambda$, to express $K$ as

$$K\, =\, 2\, \int_0^{\frac{T}{2}}\frac{r(s) \lambda(s)}{r^2(s)-1}\, ds$$
 
 When $C=C_0$, we have that $q(v_0)=0=q^\prime(v_0)$, then, we can apply the lemma 5.1 in \cite{P} and its corollary, to obtain that

$$\lim_{C\to C_0}K(C,H)=-\sqrt{2} \sqrt{1 - \frac{n H}{\sqrt{n^2 H^2-4(n-1)}}} \quad \pi=lb $$

Notice that for any $n\ge 2$ and any $H<-1$, the bound $lb<-2\pi$. By doing the substitution $v=r(s)$ in the integral  $K(C,H)$, we get that

$$K(C,H)=\int_{\tilde{t_1}(C,H)}^{\tilde{t_2}(C,H)}  \frac{2 v \, (H+(\sqrt{-C}\, v)^{-n})  }{(v^2-1)\, \sqrt{\tilde{q}(v)}}\, dv$$

When we replace $C$ by $\tilde{C}$ the integral above reduces to,

$$ K(\tilde{C},H)=\xi_n(H) $$

Using the intermediate value theorem we conclude the theorem because anytime $\xi_n(H)>-2\pi$ there exists a $C^\star\in (C_0,\tilde{C})$ such that $K(C^\star,H)=-2\pi$, therefore the map $\phi_{C^\star,H}$ is periodic in the $u$ variable, and since $C<\tilde{C}$ the function $\theta$ is injective and therefore the map $\phi_{C^\star,H}$ is an embedding.

\end{proof}

\begin{cor}
For any integer $n>1$ there exists an $H_0\le -1$ such that for any $H<H_0$ there exists an embedding with constant mean curvature $H$ from $H^{n-1}\times S^1$ to $H^{n+1}$.
\end{cor}

\begin{proof}
The corollary follows from the fact that $\lim_{H\to-\infty}\xi_n(H)=-\pi$. See item (12) in Remark (\ref{properties of q}).
\end{proof}

\begin{rem}
The integral $\xi_n$ is easy to evaluate numerically, for example

$$\xi_3(-1)=-5.97106763713693\quad  \xi_4(-1)=-4.599155062889069\quad \xi_5(-1)=-4.13016242612799$$

The following graphs suggest that for $n=3,4,5$, there exist embeddings for all $H<-1$.

\begin{figure}[h]
\centerline{\includegraphics{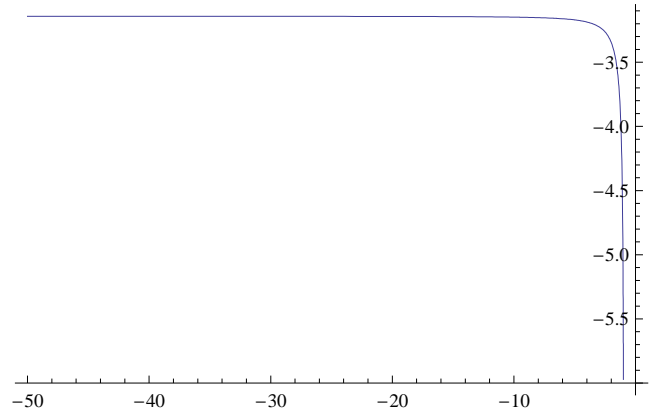}}
\caption{Graph of the function $\xi_3$ on $[-50,-1]$ }
\end{figure}
\vfill
\eject

\begin{figure}[h]
\centerline{\includegraphics{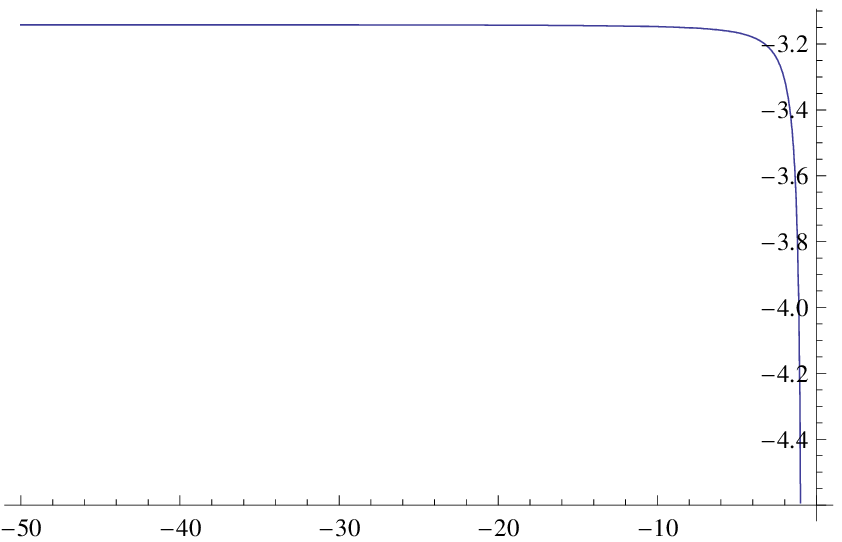}}
\caption{Graph of the function $\xi_4$ on $[-50,-1]$ }
\end{figure}

\begin{figure}[h]
\centerline{\includegraphics{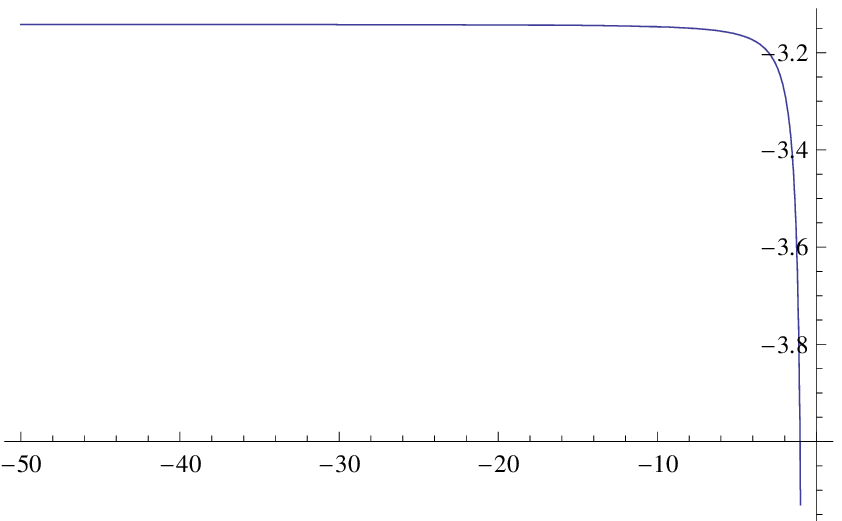}}
\caption{Graph of the function $\xi_5$ on $[-50,-1]$ }
\end{figure}
\end{rem}

\end{document}